\documentclass[secthm,seceqn,amsthm,ussrhead,12pt,reqno]{amsart}
\usepackage[utf8]{inputenc}
\usepackage[english]{babel}
\usepackage{amssymb,amsmath,amsthm,amsfonts,xcolor,enumerate,hyperref,comment,longtable,cleveref}

\usepackage{times}
\usepackage{cite}
\usepackage{pdflscape}
\usepackage{ulem}
\usepackage[mathcal]{euscript}
\usepackage{tikz}
\usepackage{hyperref}
\usepackage{cancel}
\usepackage{stmaryrd}

\usetikzlibrary{arrows}

\newtheorem{theorem}{Theorem}

\newtheorem{remark}[theorem]{Remark}

\crefrangeformat{equation}{#3(#1)#4--#5(#2)#6}

\crefname{theorem}{Theorem}{Theorems}
\crefname{lemma}{Lemma}{Lemmas}
\crefname{corollary}{Corollary}{Corollaries}
\crefname{proposition}{Proposition}{Propositions}
\crefname{definition}{Definition}{Definitions}
\crefname{example}{Example}{Examples}
\crefname{remark}{Remark}{Remarks}
\crefname{section}{Section}{Sections}
\crefname{equation}{\unskip}{\unskip}
\crefname{enumi}{\unskip}{\unskip}

\newenvironment{Proof}[1][Proof.]{\begin{trivlist}
\item[\hskip \labelsep {\bfseries #1}]}{\flushright
$\Box$\end{trivlist}}

\usepackage{stmaryrd}
\usepackage{xcolor}

\begin{document}
\noindent{\Large 
The geometric classification of nilpotent Tortkara algebras}\footnote{
The work was supported by RSF 18-71-10007. 
} \footnote{Corresponding Author: kaygorodov.ivan@gmail.com  }

   \

   {\bf  Ilya Gorshkov$^{a,b,c}$,  Ivan Kaygorodov$^{c}$ \&  Mykola Khrypchenko$^{d,e}$ \\

    \medskip
}

{\tiny

$^{a}$ Sobolev Institute of Mathematics, Novosibirsk, Russia

$^{b}$ Siberian Federal University, Krasnoyarks, Russia

$^{c}$ CMCC, Universidade Federal do ABC, Santo Andr\'e, Brazil

$^{d}$ Departamento de Matem\'atica, Universidade Federal de Santa Catarina, Florian\'opolis, Brazil

$^{e}$ Departamento de Matem\'atica, Faculdade de Ci\^{e}ncias e Tecnologia, Universidade Nova de Lisboa, Caparica, Portugal
\

\smallskip

   E-mail addresses:

\smallskip
    Ilya Gorshkov (ilygor8@gmail.com)
    
    Ivan Kaygorodov (kaygorodov.ivan@gmail.com) 
    
    Mykola Khrypchenko (nskhripchenko@gmail.com)

}

\ 

\

\ 

\noindent {\bf Abstract.}
{\it We give a geometric classification of  all $6$-dimensional nilpotent Tortkara algebras over $\mathbb C$.}

\ 

\noindent {\bf Keywords}: {\it Nilpotent algebra, Tortkara  algebra, 
Lie algebra, Malcev algebra, geometric classification, degeneration.}

\ 

\noindent {\bf MSC2010}: 17A30, 17D10.

\section*{Introduction}

An anticommutative algebra $\bf A$ is called a {\it Tortkara algebra} if it satisfies the identity 
$$
(ab)(cb)=J(a,b,c)b, \mbox{ where } J(a,b,c)=(ab)c+(bc)a+(ca)b.
$$
These algebras were introduced by Dzhumadildaev in \cite{dzhuma}. It is easy to see that each metabelian Lie algebra (i.e., $(xy)(zt)=0$) is a Tortkara algebra. Another group of examples comes from Zinbiel algebras. Recall that an algebra $\bf A$ is called a {\it Zinbiel algebra} if it satisfies the identity 
$$
(xy)z=x(yz+zy).
$$ 
Zinbiel algebras were introduced by Loday in \cite{loday} and studied in \cite{cam13, dok, dzhuma5, dzhuma19, kppv, ualbay, yau}. Under the Koszul duality, the operad of Zinbiel algebras is dual to the operad of Leibniz algebras. Zinbiel algebras are also related to Tortkara algebras \cite{dzhuma} and Tortkara triple systems \cite{brem}. More precisely, every Zinbiel algebra is a Tortkara algebra under the commutator multiplication.

There are many results related to both the algebraic and geometric 
classification
of small dimensional algebras in the varieties of Jordan, Lie, Leibniz and 
Zinbiel algebras;
for algebraic results  see, for example, \cite{ack, cfk19,  gkks, gkk,   ikm19,   ikv18, kkk18, kpv19, kv16};
for geometric results see, for example, \cite{ ack, ale, ale2, aleis, maria, bb14, BC99, cfk19, gkks, gkp, GRH, GRH2, ikm19, ikv17, ikv18, kkk18, kpv19, kppv, kpv, kv16,kv17, S90}.
Here we give a geometric classification of $6$-dimensional  nilpotent Tortkara algebras over $\mathbb C.$ Our main result is Theorem~\ref{geotort} which describes the rigid algebras in this variety.

Degenerations of algebras is an interesting subject, which has been studied in various papers.
In particular, there are many results concerning degenerations of algebras of low dimension in a  variety defined by a set of identities.
One of important problems in this direction is a description of the so-called rigid algebras. 
These algebras are of big interest, since the closures of their orbits under the action of the generalized linear group form irreducible components of the variety under consideration
(with respect to the Zariski topology). 
For example, the rigid algebras in the varieties of
all $4$-dimensional Leibniz algebras \cite{ikv17},
all  $4$-dimensional nilpotent Novikov algebras \cite{kkk18},
all  $4$-dimensional nilpotent bicommutative algebras \cite{kpv19},
all  $4$-dimensional nilpotent assosymmetric algebras \cite{ikm19},
all  $6$-dimensional nilpotent binary Lie algebras \cite{ack},
and some other has been classified.
There are fewer works in which the full information about degenerations has been found for some variety of algebras.
This problem has been solved 
for $2$-dimensional pre-Lie algebras in \cite{bb09},  
for $2$-dimensional terminal algebras in \cite{cfk19},
for $3$-dimensional Novikov algebras in \cite{bb14},  
for $3$-dimensional Jordan algebras in \cite{gkp},  
for $3$-dimensional Jordan superalgebras in \cite{maria},
for $3$-dimensional Leibniz algebras in \cite{ikv18},
for $3$-dimensional anticommutative algebras in \cite{ikv18},
for $4$-dimensional Lie algebras in \cite{BC99},
for $4$-dimensional Lie superalgebras in \cite{aleis},
for $4$-dimensional Zinbiel  algebras in \cite{kppv},
for  $3$-dimensional nilpotent algebras \cite{fkkv19},
for  $4$-dimensional nilpotent Leibniz algebras in \cite{kppv},
for  $4$-dimensional nilpotent commutative algebras \cite{fkkv19},
for  $5$-dimensional nilpotent Tortkara algebras in \cite{gkks},
for  $5$-dimensional nilpotent anticommutative algebras in \cite{fkkv19},
for  $6$-dimensional nilpotent Lie algebras in \cite{S90,GRH}, 
for  $6$-dimensional nilpotent Malcev algebras in \cite{kpv}, 
for  $7$-dimensional $2$-step nilpotent Lie algebras in \cite{ale2}, 
and for all $2$-dimensional algebras in \cite{kv16}.


\section{Definitions and notations}
Given an $n$-dimensional complex vector space ${\bf V}$, the set $\rm{Hom}({\bf V} \otimes {\bf V},{\bf V}) \cong {\bf V}^* \otimes {\bf V}^* \otimes {\bf V}$ 
is a vector space of dimension $n^3$. This space has a structure of the affine variety $\mathbb{C}^{n^3}.$ Fix a basis $e_1,\dots,e_n$ of ${\bf V}$. Every $\mu\in \rm{Hom}({\bf V} \otimes {\bf V},{\bf V})$ is determined by the $n^3$ structure constants $c_{i,j}^k\in\mathbb{C}$ such that
$\mu(e_i\otimes e_j)=\sum\limits_{k=1}^nc_{i,j}^ke_k$. A subset of $\rm{Hom}({\bf V} \otimes {\bf V},{\bf V})$ is {\it Zariski-closed} if it can be defined by a set of polynomial equations in the variables $c_{i,j}^k$ ($1\le i,j,k\le n$).

Let $T$ be a set of polynomial identities.
All algebra structures on ${\bf V}$ satisfying polynomial identities from $T$ form a Zariski-closed subset of the variety $\rm{Hom}({\bf V} \otimes {\bf V},{\bf V})$. We denote this subset by $\mathbb{L}(T)$.
The general linear group $\rm{GL}({\bf V})$ acts on $\mathbb{L}(T)$ by conjugation:
$$ (g * \mu )(x\otimes y) = g\mu(g^{-1}x\otimes g^{-1}y)$$ 
for $x,y\in {\bf V}$, $\mu\in \mathbb{L}(T)\subset \rm{Hom}({\bf V} \otimes {\bf V},{\bf V})$ and $g\in \rm{GL}({\bf V})$.
Thus, $\mathbb{L}(T)$ is decomposed into $\rm{GL}({\bf V})$-orbits that correspond to the isomorphism classes of algebras. 
Let $O(\mu)$ denote the orbit of $\mu\in\mathbb{L}(T)$ under the action of $\rm{GL}({\bf V})$ and let $\overline{O(\mu)}$ denote the Zariski closure of $O(\mu)$.

Let ${\bf A}$ and ${\bf B}$ be two $n$-dimensional algebras satisfying identities from $T$ and $\mu,\lambda \in \mathbb{L}(T)$ represent ${\bf A}$ and ${\bf B}$ respectively.
We say that ${\bf A}$ {\it degenerates} to ${\bf B},$ and write ${\bf A}\to {\bf B},$ if $\lambda\in\overline{O(\mu)}$.
In this case   $\overline{O(\lambda)}\subset\overline{O(\mu)}$. Hence, the definition of a degeneration does not depend on the choice of $\mu$ or $\lambda$. If ${\bf A}\not\cong {\bf B}$, then ${\bf A}\to {\bf B}$ is a {\it proper degeneration}. We write ${\bf A}\not\to {\bf B}$ if $\lambda\not\in\overline{O(\mu)}$. Note also that if ${\bf A}\to {\bf B}$ and ${\bf B}\to {\bf C}$, then ${\bf A}\to {\bf C}$.

Let ${\bf A}$ be represented by $\mu\in\mathbb{L}(T)$. Then  ${\bf A}$ is said to be {\it rigid} in $\mathbb{L}(T)$ if $O(\mu)$ is an open subset of $\mathbb{L}(T)$.
 Recall that a subset of a variety is  {\it irreducible} if it cannot be represented as a union of two non-trivial closed subsets. 
 A maximal irreducible closed subset of a variety is  an {\it irreducible component}.
 It is well known that every affine variety can be represented as a finite union of its irreducible components in a unique way.
Observe that ${\bf A}$ is rigid in $\mathbb{L}(T)$ if and only if $\overline{O(\mu)}$ is an irreducible component of $\mathbb{L}(T)$.





We use the methods applied to Lie algebras in \cite{GRH2}.
If ${\bf A}\to {\bf B}$ and ${\bf A}\not\cong {\bf B}$, then $\dim \mathfrak{Der}({\bf A})<\dim \mathfrak{Der}({\bf B})$, where $\mathfrak{Der}({\bf A})$ is the Lie algebra of derivations of ${\bf A}$. We will compute the dimensions of algebras of derivations and will check the assertion ${\bf A}\to {\bf B}$ only for such ${\bf A}$ and ${\bf B}$ that $\dim \mathfrak{Der}({\bf A})<\dim \mathfrak{Der}({\bf B})$. 


To prove degenerations, we  construct families of matrices parametrized by $t$. Namely, let ${\bf A}$ and ${\bf B}$ be two algebras represented by the structures $\mu$ and $\lambda$ from $\mathbb{L}(T)$ respectively. Let $e_1,\dots, e_n$ be a basis of $\bf  A$ and let $c_{i,j}^k$ ($1\le i,j,k\le n$) be the structure constants of $\lambda$ in this basis. If there exist $a_i^j(t)\in\mathbb{C}$ ($1\le i,j\le n$, $t\in\mathbb{C}^*$) such that $E_i^t=\sum\limits_{j=1}^na_i^j(t)e_j$ ($1\le i\le n$) form a basis of ${\bf V}$ for every $t\in\mathbb{C}^*$, and the structure constants of $\mu$ in the basis $E_1^t,\dots, E_n^t$ are  functions $c_{i,j}^k(t)$ such that $\lim\limits_{t\to 0}c_{i,j}^k(t)=c_{i,j}^k$, then ${\bf A}\to {\bf B}$. In this case  $E_1^t,\dots, E_n^t$ is called a {\it parametrized basis} for ${\bf A}\to {\bf B}$.

\section{The geometric classification of $6$-dimensional nilpotent Tortkara  algebras}
The geometric classification of $6$-dimensional nilpotent Tortkara algebras is  based on the description of all degenerations of $6$-dimensional nilpotent Malcev algebras~\cite{kpv}.
Note that the algebras $g_6$ and $g_8$ from~\cite{kpv} are not Tortkara.
Hence, every  $6$-dimensional nilpotent Tortkara-Malcev algebra  degenerates from one of the following algebras:

$$
\begin{array}{llllllllllllllllllllllllll}

g_5 &:&  e_1e_2 = e_3,& e_1e_3 = e_4,& e_1e_4 = e_5,& e_1e_5 = e_6,& e_2e_3 = e_5,& e_2e_4 = e_6,\\   
    
{\bf M}_{6}^{\epsilon}      &:&  e_1e_2=e_3, &e_1e_3=e_5, &e_1e_5=e_6,&e_2e_4=\epsilon e_5, &e_3e_4=e_6.  
\end{array}
$$
It remains to take into account all $6$-dimensional nilpotent  non-Malcev Tortkara algebras. These algebras were described in~\cite{gkk}. We recall in the next theorem their classification.

\begin{theorem}
Let $\mathbb T$ be a $6$-dimensional nilpotent non-Malcev Tortkara algebra over $\mathbb C$. 
Then $\mathbb T$ is isomorphic to one of the following algebras:

$$
\begin{array}{lllllllllllllllllllllllll}

{\mathbb T}_{00}^6 & : &  e_1e_2=e_3, &  e_1e_3=e_4, & e_2e_4=e_5;\\

{\mathbb T}_{01}^6 & : &  
e_1e_2=e_3,& e_1e_3=e_4, & e_1e_4=e_5, & e_2e_3=e_5, & e_2e_4=e_6;\\

{\mathbb T}_{02}^6 & : & 
e_1e_2=e_3,& e_1e_3=e_4,& e_2e_3=e_5,& e_2e_4=e_6;\\

{\mathbb T}_{03}^6 & : &    
e_1e_2=e_3, & e_1e_3=e_4, & e_1e_4=e_5, & e_2e_4=e_6;\\

{\mathbb T}_{04}^6 & : &    
e_1e_2=e_3, & e_1e_3=e_4,& e_1e_5=e_6,& e_2e_4=e_6;\\ 

{\mathbb T}_{05}^6 & : &
e_1e_2=e_3, & e_1e_3=e_4 , & e_2e_3=e_6, & e_4e_5=e_6;\\

{\mathbb T}_{06}^6 & : &
e_1e_2=e_3, & e_1e_3=e_4, & e_2e_4=e_6, & e_3e_5=e_6;\\

{\mathbb T}_{07}^6 & : &
e_1e_2=e_3, & e_1e_3=e_4, &  e_4e_5=e_6;\\

{\mathbb T}_{08}^6 & : &  
e_1e_2=e_3, & e_1e_3=e_4, & e_1e_4=e_6, & e_1e_5=-e_6, & e_2e_3=e_5, & e_2e_4=e_6;\\

{\mathbb T}_{09}^6(\alpha) & : &  
e_1e_2=e_3, & e_1e_3=e_4, & e_1e_5=(\alpha+1) e_6, & e_2e_3=e_5, & e_2e_4=\alpha e_6;\\

{\mathbb T}_{10}^6 & : &  
e_1e_2=e_3, & e_1e_3=e_6, & e_1e_4=e_5, & e_2e_3=e_5, & e_4e_5=e_6;\\

{\mathbb T}_{11}^6 & : & 
e_1e_2=e_3, & e_1e_4=e_5, & e_1e_5=e_6, &  e_2e_3=e_5;\\

{\mathbb T}_{12}^6 & : & 
e_1e_2=e_3, & e_1e_4=e_5, & e_1e_5=e_6, &  e_2e_3=e_5, &  e_2e_4=e_6;\\

{\mathbb T}_{13}^6 & : & 
e_1e_2=e_3, & e_1e_4=e_5, & e_1e_5=e_6, &  e_2e_3=e_5, & e_3e_4=e_6;\\

{\mathbb T}_{14}^6 & : & 
e_1e_2=e_3, & e_1e_4=e_5, & e_2e_3=e_5, & e_4e_5=e_6;\\

{\mathbb T}_{15}^6 & : & 
e_1e_2=e_3, & e_1e_3=e_4, & e_1e_4=e_5, & e_1e_5=e_6, & e_2e_4=e_6;\\

{\mathbb T}_{16}^6 & : &  
e_1e_2=e_3, & e_1e_3=e_4, & e_1e_4=e_5, & e_2e_5=e_6;\\

\mathbb{T}_{17}^6 & : &
e_1e_2=e_3, & e_1e_3=e_4, & e_1e_4=e_5, & e_2e_3=e_5, & e_2e_5=e_6;\\

\mathbb{T}_{18}^6(\alpha) & : &
e_1e_2=e_3, & e_1e_3=e_4, & e_1e_4=e_5, & e_1e_5=(\alpha+1)e_6, & e_2e_3=e_5, & e_2e_4=\alpha e_6;\\

{\mathbb T}_{19}^6 & : & 
e_1e_2=e_3, & e_1e_3=e_4, & e_1e_5=e_6, & e_2e_4=e_5, & e_3e_4=e_6.

\end{array}
$$

All listed algebras are non-isomorphic, except  
 ${\mathbb T}^6_{09}(\alpha) \cong {\mathbb T}^6_{09}(-\alpha-1).$

\end{theorem}

\begin{remark}
There is only one $6$-dimensional nilpotent non-metabelian Tortkara algebra: ${\mathbb T}_{19}^6$. 
Obviously, ${\mathbb T}_{19}^6$ is rigid.
\end{remark}

The main result of the present paper is the following theorem.

\begin{theorem}\label{geotort}
The variety of $6$-dimensional nilpotent Tortkara algebras over $\mathbb C$ has $3$ irreducible components
defined by the rigid algebras ${\mathbb T}_{10}^6,$ ${\mathbb T}_{17}^6$ and ${\mathbb T}_{19}^6.$
\end{theorem}

\begin{Proof}
The algebras ${\mathbb T}_{10}^6,$ ${\mathbb T}_{17}^6$ and ${\mathbb T}_{19}^6$ are rigid, because they have the minimal dimension of derivation algebras among all $6$-dimensional nilpotent Tortkara algebras, which is $7$. The rest of $6$-dimensional nilpotent Tortkara algebras degenerate from (at least) one of these algebras, as shown in the table below.

\end{Proof}


{\tiny 
$$
\begin{array}{|lcl|lll|}

\multicolumn{6}{c}{ \mbox{ {\bf Table.}
\mbox{\it  Degenerations of $6$-dimensional nilpotent Tortkara algebras.}}} \\
\multicolumn{6}{c}{}\\

\hline

{\mathbb  T}^6_{19}    &\to& g_{5} & 
E_1^t = t e_1 + t^{-2} e_2 - t^{-3} e_3 -  t^{-4} e_4, & E_2^t = e_2 + t^{-1} e_3 + t^{-2} e_4, & E_3^t = t e_3 + e_4 + 2 t^{-4} e_5 \\ 
&&& E_4^t = t^2 e_4 + t^{-2} e_5 + 2 t^{-3}e_6, & E_5^t = e_5, & E_6^t = te_6 \\
\hline

{\mathbb  T}^6_{19}    &\to& {\bf M}_{6}^{\epsilon}    & 
E_1^t = t e_1 - t^{-1} e_2 - t^{-2}e_3, & E_2^t = e_2 - \frac{ \epsilon-1}{t}e_3 , & E_3^t = e_4 + \epsilon t^{-3}e_5 \\ 
&&& E_4^t = e_3 - \epsilon t^{-1}e_4 - \epsilon t^{-4}e_5, & E_5^t = -t^{-1} e_5 + \frac{\epsilon-1}{ t^2} e_6, & E_6^t = -e_6\\
\hline

{\mathbb  T}^6_{19}    &\to& {\mathbb T}^6_{00}    & 
E_1^t = e_1, & E_2^t = e_2, & E_3^t = e_3 \\
&&& E_4^t = e_4, & E_5^t = e_5, & E_6^t = t^{-1}e_6\\
\hline

{\mathbb  T}^6_{19}    &\to& {\mathbb T}^6_{01}    & 
E_1^t = e_1+e_2-t^{-2}e_3, & E_2^t = t e_2+e_3+t^{-2} e_4, & E_3^t = t e_3+e_4+t^{-2} e_5-t^{-4} e_6 \\
&&& E_4^t = t e_4+e_5, & E_5^t = t e_5+(1 - t^{-1})e_6, & E_6^t = e_6\\
\hline

{\mathbb  T}^6_{19}    &\to& {\mathbb T}^6_{02}    & 
E_1^t = te_1, & E_2^t = e_2+t^{-1} e_3+t^{-3}e_4, & E_3^t = te_3+e_4 \\
&&& E_4^t = t^2 e_4, & E_5^t = e_5+(t^{-1} - t^{-2})e_6, & E_6^t = e_6 \\
\hline

{\mathbb  T}^6_{19}    &\to& {\mathbb T}^6_{03}    & 

E_1^t = t e_1 + t^{-2}e_2 - t^{-4}e_3, & E_2^t = e_2, & E_3^t = te_3\\
&&& E_4^t = t^2 e_4, & E_5^t = e_5 - t^{-2}e_6, & E_6 = e_6\\
\hline

{\mathbb  T}^6_{19}    &\to& {\mathbb T}^6_{04}    & 
E_1^t = t e_1, & E_2^t = e_2, & E_3^t = t e_3 \\
&&& E_4^t = t^2 e_4, & E_5^t = e_5 - t^{-1}e_6, & E_6 = t e_6\\
\hline

{\mathbb  T}^6_{10}    &\to& {\mathbb T}^6_{05}    & 
E_1^t = e_1 + t^{-1}e_2 - te_5, & E_2^t = e_2, &  E_3^t = e_3\\ 
&&&E_4^t = t^{-1}e_5 + e_6, & E_5^t = t^2 e_4, & E_6^t = -t e_6\\
\hline

{\mathbb  T}^6_{05}    &\to& {\mathbb T}^6_{06}    & 
E_1^t = -e_1, & E_2^t = e_2 + t^{-1}e_3 + te_5, & E_3^t = -e_3 - t^{-1}e_4, \\ 
&&& E_4^t = e_4, & E_5^t = t^2 e_5, & E_6^t = -t e_6\\
\hline

{\mathbb  T}^6_{05}    &\to& {\mathbb T}^6_{07}    & 
 E_1^t =  te_1, &   E_2^t =  e_2, & E_3^t = t e_3 \\
&&& E_4^t = t^{2} e_4, & E_5^t = t^{-2} e_5, & E_6^t = e_6\\
\hline

{\mathbb  T}^6_{19}    &\to& {\mathbb T}^6_{08}    & 
E_1^t = t e_1 - t^{-1}e_3, & E_2^t = e_2 + t^{-1}e_3 - t^{-2}e_4, &  E_3^t = t e_3 + e_4 +  t^{-3}e_6, \\
&&& E_4^t = t^2 e_4 - t^{-1}e_6, & E_5^t = e_5 + 2t^{-1}e_6, & E_6 = -t e_6 \\
\hline

{\mathbb  T}^6_{19}    &\to& {\mathbb T}^6_{09}(\alpha)    & 
E_1^t = te_1, & E_2^t = e_2 + t^{-1} e_3 + \frac{\alpha}{(\alpha + 1)t^2} e_4, & E_3^t = te_3 + e_4, \\ 
&&& E_4^t = t^2 e_4, & E_5^t = e_5 + \frac{1}{(\alpha + 1)t}e_6, & E_6^t = \frac{t}{\alpha+1} e_6 \\
\hline

{\mathbb  T}^6_{12}    &\to& {\mathbb T}^6_{11}    & 
E_1^t = e_1, & E_2^t = t e_2, & E_3^t = t  e_3 \\
&&& E_4^t = t^{2} e_4, & E_5^t = t^{2} e_5, & E_6^t = t^{2}e_6\\
\hline

{\mathbb  T}^6_{19}    &\to& {\mathbb T}^6_{12}    & 
E_1^t = e_1 + e_2 - t^{-1}e_3, & E_2^t = t e_2 + t^{-1}e_3 +  t^{-2}e_4, & E_3^t = te_3 + t^{-1}e_4 + t^{-2}e_5 -  t^{-3}e_6 \\
&&& E_4^t = e_4 +  t^{-2}e_5, & E_5^t = e_5 + (t^{-2}-t^{-1})e_6, & E_6 = e_6\\
\hline

{\mathbb  T}^6_{10}    &\to& {\mathbb T}^6_{13}    & 
E_1^t = e_1 + t^{-1}e_3 + t e_4, & E_2^t = t^2 e_2, & E_3^t = t^2 e_3 - te_5 \\ 
&&& E_4^t = t^4 e_4, & E_5^t = t^4 e_5, & E_6^t = t^5 e_6\\
\hline

{\mathbb  T}^6_{10}    &\to& {\mathbb T}^6_{14}    & 
E_1^t = te_1,& E_2^t = e_2,& E_3^t = t e_3 \\
&&& E_4^t = e_4, & E_5^t = t e_5, & E_6^t = t e_6 \\
\hline

{\mathbb  T}^6_{19}    &\to& {\mathbb T}^6_{15}    & 
E_1^t = te_1 + t^{-2} e_2 - t^{-3}e_3, & E_2^t = e_2, & E_3^t = t e_3 \\
&&& E_4^t = t^2 e_4, & E_5^t = e_5 - t^{-1} e_6, & E_6^t = te_6 \\
\hline

{\mathbb  T}^6_{17}    &\to& {\mathbb T}^6_{16}    & 
E_1^t = t^{-1}e_1, & E_2^t =  e_2, & E_3^t = t^{-1} e_3 \\
&&& E_4^t = t^{-2} e_4, & E_5^t = t^{-3} e_5, &  E_6^t = t^{-3}e_6\\
\hline

{\mathbb  T}^6_{19} &\to & {\mathbb  T}^6_{18}(\alpha) &

E_1^t = te_1 + t^{-2} e_2 - \frac{\alpha}{(\alpha+1)t^3} e_3,&
E_2^t = e_2 + t^{-1} e_3 + \frac{\alpha}{ (\alpha + 1)t^2}  e_4,&
E_3^t = t e_3 + e_4 + \frac{\alpha}{(\alpha + 1)t^4}e_5 - \frac{\alpha^2}{(\alpha+1)^2t^5} e_6 \\
&&& 
E_4^t= t^2e_4 + t^{-2}e_5, &   E_5^t = e_5 + \frac{1}{(\alpha + 1)t}e_6, & E_6^t = \frac{t}{\alpha+1} e_6\\ \hline

\end{array}$$
}

\end{document}